\documentclass[12pt]{amsart}
\usepackage{amscd}

\DeclareFontEncoding{OT2}{}{} 
\newcommand{\textcyr}[1]{%
 {\fontencoding{OT2}\fontfamily{cmr}\fontseries{m}\fontshape{n}\selectfont #1}}

\newcommand{\Sha}{{\mbox{\textcyr{Sh}}}}
\newcommand{\PIC}{{\operatorname{\bf Pic}}}

\newcommand{\real}{{\operatorname{real}}}
\newcommand{\tbar}{{\bar{t}}}

\newcommand{\Aff}{{\mathbf A}}
\newcommand{\A}{{\mathcal A}}
\newcommand{\B}{{\mathcal B}}
\newcommand{\C}{{\mathcal C}}

\newcommand{\Ebar}{{\overline{E}}}
\newcommand{\E}{{\mathcal E}}

\newcommand{\V}{{\mathcal V}}
\newcommand{\X}{{\mathcal X}}
\newcommand{\Y}{{\mathcal Y}}
\newcommand{\Ghat}{{\hat{G}}}
\newcommand{\Atilde}{{\tilde{\A}}}
\newcommand{\Xtilde}{{\tilde{\X}}}
\newcommand{\Ctilde}{{\tilde{C}}}
\newcommand{\OO}{{\mathcal O}}

\newcommand{\mm}{{\mathfrak m}}
\newcommand{\Xbar}{\overline{\X}}

\newcommand{\Q}{{\mathbf Q}}
\newcommand{\Qbar}{\overline{\Q}}
\newcommand{\kbar}{{\overline{k}}}

\newcommand{\Z}{{\mathbf Z}}
\newcommand{\R}{{\mathbf R}}
\newcommand{\PP}{{\mathbf P}}
\newcommand{\PPdual}{{\check{\PP}}}

\newcommand{\F}{{\mathbf F}}
\newcommand{\G}{{\mathbf G}}

\newcommand{\Ak}{{{\mathbb A}_k}}
\newcommand{\Fbar}{{\overline{\F}}}

\newcommand{\et}{{\operatorname{\acute{e}t}}}
\newcommand{\HH}{\operatorname{H}}
\newcommand{\RR}{\operatorname{R}}
\newcommand{\ev}{\operatorname{ev}}

\newcommand{\Hom}{\operatorname{Hom}}

\newcommand{\Gal}{\operatorname{Gal}}

\newcommand{\Br}{\operatorname{Br}}

\newcommand{\Div}{\operatorname{Div}}
\newcommand{\Pic}{\operatorname{Pic}}
\newcommand{\Res}{\operatorname{Res}}
\newcommand{\Spec}{\operatorname{Spec}}

\newcommand{\nichts}{{\left.\right.}}
\newcommand{\isom}{\cong}

\newtheorem{theorem}{Theorem}
\newtheorem{lemma}[theorem]{Lemma}
\newtheorem{cor}[theorem]{Corollary}
\newtheorem{prop}[theorem]{Proposition}

\theoremstyle{definition}

\theoremstyle{remark}
\newtheorem{rem}{Remark$\!\!$}	
\newtheorem{rems}{Remarks$\!\!$}	

\begin{document}

\title[Families of Shafarevich-Tate group elements]
{Algebraic families of nonzero elements of Shafarevich-Tate groups}
\subjclass{Primary 11G10; Secondary 11G30, 11G35, 14H40, 14J27}
\keywords{Shafarevich-Tate group, Brauer-Manin obstruction, Hasse principle,
cubic surface, Cassels-Tate pairing, Lefschetz pencil}
\author{Jean-Louis Colliot-Th\'el\`ene}
\address{C.N.R.S., Math\'ematiques, B\^atiment 425, Universit\'e de Paris-Sud,
F-91405 Orsay, France}
\email{colliot@math.u-psud.fr}
\author{Bjorn Poonen}
\thanks{Most of the research for this paper was done while the authors
were both enjoying the hospitality of the Isaac Newton Institute,
Cambridge, England.
The first author is a researcher at C.N.R.S.
The second author is partially supported by NSF grant DMS-9801104,
a Sloan Fellowship, and a Packard Fellowship.}
\address{Department of Mathematics, University of California,
Berkeley, CA 94720-3840, USA}
\email{poonen@math.berkeley.edu}
\date{January 1, 1999}


\maketitle

\section{Introduction}

The purpose of this note is to show that there exist
non-trivial families of algebraic varieties for which
all the fibers above rational points
(or even above points of odd degree)
are torsors of abelian varieties
representing nonzero elements of their Shafarevich-Tate groups.
More precisely, we will prove the following theorem.
(Throughout this paper, $\PP^1$ unadorned denotes $\PP^1_\Q$.)

\begin{theorem}
\label{main}
There exists
an open subscheme $U$ of $\PP^1$ containing
all closed points of odd degree,
and there exist
smooth projective geometrically integral varieties $\A$ and $\X$ over
$\Q$ equipped with dominant morphisms $\pi_\A$ and $\pi_\X$ to $\PP^1$
such that
\begin{itemize}
	\item[(1)] $\A_U:=\pi_\A^{-1}(U)$ is an abelian scheme
			of relative dimension 2 over $U$.
	\item[(2)] $\X_U:=\pi_\X^{-1}(U)$ is an $\A_U$-torsor over $U$.
	\item[(3)] $\pi_\X(\X(k_v))=\PP^1(k_v)$ for every local field
		$k_v \supset \Q$ (archimedean or not).
	\item[(4)] $\X$ has no zero-cycle of odd degree over $\Q$.
	\item[(5)] $\A$ is a non-constant family (i.e., there exist
			$u,v \in U(\Qbar)$ such that the fibers
			$\A_u$ and $\A_v$ are not isomorphic
			over $\Qbar$)
	\item[(6)] The generic fiber of $\A \rightarrow \PP^1$
		is absolutely simple (i.e., simple over $\overline{\Q(t)}$).
\end{itemize}
\end{theorem}

\begin{rems}
$\nichts$
\begin{enumerate}
	\item Conditions (1) through (4) imply that
		for each closed point $s \in U$ of odd degree
		(in particular, for all $s \in \PP^1(\Q)$),
		the fiber $\X_s$ represents a nonzero element of
		the Shafarevich-Tate group $\Sha(\A_s)$.
	\item Let $E$ be an elliptic curve over $\Q$,
		and suppose that $C$ is a genus~1 curve over $\Q$
		representing an element of $\Sha(E)$ whose order is even.
		Let $\B$ be a smooth projective geometrically integral
		variety over $\Q$ equipped
		with a morphism $\pi_\B: \B \rightarrow \PP^1$
		making $\pi_\B^{-1}(U)$ a
		non-constant abelian scheme over $U$,
		$U$ being as above.
		Then $\A:=E \times_\Q \B$ and $\X:=C \times_\Q \B$,
		equipped with the morphisms to $\PP^1$ obtained by composing
		the second projection with $\B \rightarrow \PP^1$,
		satisfy (1) through (5).
		Condition (6) is designed to rule out such ``trivial''
		examples.
\end{enumerate}
\end{rems}

The varieties $\A$ and $\X$ promised by the theorem will be
smooth compactifications (over $\Q$) of $\PIC^0_{\C/\PP^1}$ and
$\PIC^1_{\C/\PP^1}$
for a relative curve $\C \rightarrow \PP^1$ of genus~2.
(For the theory of the relative Picard functor, we refer the reader to
\cite{blr}.)
The key ingredient
which
will let us prove the non-existence
of rational points of $\X$ over fields of odd degree (condition~(4))
is the formula from~\cite{poonenstoll}
which
gives the value of Cassels-Tate pairing when the
torsor $\PIC^{g-1}$ of $\PIC^0$ for a genus $g$ curve
over a global field is paired with itself.
Using that same formula,
we will show
that the non-existence of zero-cycles of degree one on $\X$
can be explained by a Brauer-Manin obstruction.

We will also give an example of relative dimension~1,
but only for fibers above rational points,
not above all points of odd degree.

\begin{theorem}
\label{genus1}
There exists an open subscheme $U$ of $\PP^1$ containing $\PP^1(\Q)$,
and smooth projective geometrically integral varieties $\E$ and $\X$ over $\Q$
equipped with dominant morphisms $\pi_\E$ and $\pi_\X$ to $\PP^1$
such that
\begin{itemize}
	\item[(1)] $\X_U \rightarrow U$ is a proper smooth family whose fibers
		are geometrically integral curves of genus~1.
	\item[(2)] $\E_U = \PIC^0_{\X_U/U}$,
		so $\E$ is an elliptic surface over $\PP^1$, smooth above $U$.
	\item[(3)] $\pi_\X(\X(\Q_p))=\PP^1(\Q_p)$ for all $p \le \infty$.
	\item[(4)] $\X(\Q)$ is empty.
	\item[(5)] The $j$-invariant of $\E_U \rightarrow U$ is a non-constant
		function on $U$.
\end{itemize}
\end{theorem}

The proof of this begins by taking
a Lefschetz pencil in a cubic surface violating the Hasse principle.
A suitable base change of this will give us $\X \rightarrow \PP^1$.
When the failure of the Hasse principle for the cubic surface
is due to a Brauer-Manin obstruction\footnote{
   It has been conjectured that the Brauer-Manin obstruction
   is the only obstruction to the Hasse principle for cubic surfaces.
   In any case, there are examples of cubic surfaces violating
   the Hasse principle because of such an obstruction.
   See~\ref{review} and~\ref{violation} below.},
we can construct a second family $\Y_U \rightarrow U$ of
genus~1 curves such that
$\PIC^0_{\Y_U/U} = \E_U$,
the analogue of~(3) holds,
and the Cassels-Tate pairing satisfies $\langle \X_t,\Y_t \rangle=1/3$ for
all $t \in U(\Q)$.
While proving the last statement, we are led to prove an auxiliary result
(Lemma~\ref{oneelement}) which may be of independent interest:
if there is a Brauer-Manin obstruction
to the Hasse principle for a smooth cubic surface $V$,
then only one element $\alpha \in \Br(V)$
is needed to create the obstruction.

 Throughout this paper,  \'etale cohomology of a scheme $X$
 with values in a sheaf $\mathcal F$ is denoted $\HH^{*}_\et(X,{\mathcal F})$.
 The notation for the Galois cohomology of (the spectrum of) a
 field $k$ will simply be $\HH^{*}(k,{\mathcal F})$.
The (cohomological) Brauer group $\Br(X)$ of a scheme $X$ is the group
$\HH^2_\et(X,\G_m)$.
The Shafarevich-Tate group $\Sha(A)$
of an abelian variety $A$ over a number field $k$
is defined as the kernel of the natural map
$\HH^1(k,A) \rightarrow \prod_v \HH^1(k_v,A)$,
where the product is over all places $v$ of $k$,
and $k_v$ denotes the completion.
If $k$ is a field, then $\kbar$ denotes an algebraic closure of $k$,
and if $X$ is a $k$-scheme, then $\overline{X} := X \times_k \kbar$.
Given an integral (i.e., reduced and irreducible) $k$-variety, we write
$k(X)$ for the function field of $X$. If $X/k$ is geometrically
integral, we denote by $\kbar(X)$ the function field
of ${\overline X}$.

\section{A family of abelian surfaces}
\label{abeliansurfaces}

In this section we prove Theorem~\ref{main}.

\subsection{Minimal models of hyperelliptic curves}

We need a simple lemma describing the minimal model
of certain hyperelliptic curves.
Recall~\cite{delignemumford} that a {\em stable curve} over a base $S$
is a proper flat morphism $\C \rightarrow S$ whose geometric fibers
are reduced, connected, $1$-dimensional schemes $\C_s$ such that
$\C_s$ has only ordinary double points,
and any nonsingular rational component of $\C_s$ meets
the other components in more than~2 points.

\begin{lemma}
\label{semistable}
Let $R$ be a Dedekind domain with $2 \in R^\ast$,
and let $K$ be its fraction field.
Let $f(x,z) \in R[x,z]$ be homogeneous of even degree $2n \ge 4$,
and assume that its discriminant generates a squarefree ideal of $R$.
Let $C/K$ be the smooth projective model of the affine
curve defined by $y^2=f(x,1)$ over $K$.
Then the minimal proper regular model $\C/R$ of $C$ is stable
and can be obtained
by glueing the affine curves $y^2=f(x,1)$ and $Y^2=f(1,z)$ over $R$
along the open subsets where $x \not=0$ and $z \not=0$ respectively,
using the identifications $z=x^{-1}$ and $Y=x^{-n}y$.
In particular, its geometric fibers are integral.
\end{lemma}

\begin{proof}
The model $\C$ described is finite over $\PP^1_R$
(which manifests itself here
as two copies of $\Aff^1_R$ with variables $x$, $z$,
glued along $x \not= 0$ and $z \not=0$ using the identification $z=x^{-1}$),
so $\C$ is proper over $R$.
The regularity and stability follow from Corollaire~6
and Remarque~18 in~\cite{liu}.
The integrality of the fibers follows from the explicit construction,
and implies that $\C$ is the minimal model.
\end{proof}

\subsection{The relative curve}
\label{relativecurve}

Let $C$ be the smooth projective model of the curve
\begin{equation}
\label{model1}
	y^2 = - \left[ x^6 + x^5 + t^5 x + (8t^6+7) \right]
\end{equation}
over $\Q(t)$.
Let $f_t(x) \in \Q[t][x]$ denote the right hand side of~(\ref{model1}).
PARI shows that the discriminant $\Delta(t)$ of $f_t(x)$ with respect to $x$
is an irreducible polynomial of degree 30 in $\Q[t]$.
In particular $\Delta(t) \not\equiv 0$, so $C$ is of genus~2.
Let $s_0 \in \Aff^1_\Q \subset \PP^1$ be the
closed point corresponding to $\Delta(t)$,
and let $U = \PP^1 \setminus \{s_0\}$.

Since $\Delta(t)$ is squarefree, the minimal model $\C^0$ over $\Spec \Q[t]$
is described by Lemma~\ref{semistable}.
Let $T=t^{-1}$.
Then $C$ is birational over $\Q(t)$ to the affine curve
\begin{equation}
\label{model2}
	Y^2 = - \left[ X^6 + T X^5 + X + (8+7 T^6) \right],
\end{equation}
via the change of variables $X=x/t$ and $Y=y/t^3$.
The minimal model $\C^1$ over $\Spec \Q[T]$
is given by Lemma~\ref{semistable} again.
Glueing $\C^0$ and $\C^1$ gives the minimal model $\C \rightarrow \PP^1$.
Then $\C \rightarrow \PP^1$ is stable, and has good reduction over $U$,
since this can be checked locally on the base.
Moreover $\C$ is smooth and proper over $\Q$.

\subsection{Local points on the fibers}

In what follows, we will specialize $\C$ at points in $\PP^1(k)$
where $k$ is a field containing $\Q$.
Although $t$ above denoted an indeterminate,
we will abuse notation below
by letting $t$ denote also the specialized value in $k \cup \{\infty\}$,
and the corresponding $k$-valued point $\Spec k \rightarrow \PP^1$.
If $t \not= \infty$
(resp.~$t \not= 0$), then the fiber $\C_t$ is
birational to the curve over $k$ given by the equation~(\ref{model1})
(resp.~(\ref{model2})).

\begin{lemma}
\label{localdivisors}
$\nichts$
\begin{enumerate}
	\item If $t \in U(\R)$,
		then $\C_t$ has no (real) divisor of degree~1.
	\item If $k$ is a finite extension of $\Q_p$ for some finite
		prime $p$, and if $t \in U(k)$,
		then $\C_t(k) \not= \emptyset$.
\end{enumerate}
\end{lemma}

\begin{proof}
For $t,x \in \R$, the weighted AM-GM inequality gives
	$$\frac{1}{6}x^6 + \frac{5}{6}t^6 \ge |t^5 x| \ge -t^5 x
	\qquad \text{and} \qquad
	\frac{5}{6}x^6 + \frac{1}{6} \ge |x^5| \ge -x^5.$$
Combining these shows that
	$$ x^6 + x^5 + t^5 x + (8 t^6+7) > 0. $$
Thus $\C$ has a dense open subscheme with no real points.
But $\C\times_{\Q}\R$
is smooth over $\R$, so $\C(\R)=\emptyset$.
Hence $\C$ and $\C_t$ for $t \in U(\R)$
have no real zero-cycles of odd degree.

{}From now on, we assume $t \in U(k)$ where $[k:\Q_p]<\infty$,
and we let $v_p$ be the $p$-adic valuation on $k$,
normalized so that $v_p(p)=1$.
By convention, if $t=\infty$, then $v_p(t)=-\infty<0$.
First suppose $p=2$.
If $v_2(t) < 0$, then Hensel's Lemma shows that~(\ref{model2})
has a $k$-point with $Y=0$, and with $X$ near~0.
If $v_2(t) \ge 0$, then~(\ref{model1}) has a $k$-point with $x=0$.

{}From now on, we assume $p$ is odd.
Let $\F_q$ be the residue field of $k$.
Suppose $v_p(t) \ge 0$.
Then let $f(x) \in \F_q[x]$ denote the reduction of the right hand side
of~(\ref{model1}).
If $f(x)$ is a square in $\Fbar_p[x]$,
then equating coefficients of $x^6, x^5, x^4, x^3$ shows that
the reduction would have to be $-(x^3+x^2/2-x/8+1/16)^2$,
but equating the coefficients of $x^2,x^1,x^0$ gives the
inconsistent system
	$$0 = -5/64, \qquad \tbar^5 = -1/64, \qquad 7+8\tbar^6 = 1/256,$$
of equations in $\Fbar_p$, where $\tbar$ denotes the reduction of $t$.
Hence $f(x)$ is not the square of any polynomial in $\Fbar_p[X]$.
If instead $v_p(t) < 0$, then let $f(x)=-(X^6+X+8) \in \F_q[x]$
be the reduction of the right hand side of~(\ref{model2}),
which again is not a square in $\Fbar_p[X]$.

In either case (no matter what the value of $v_p(t)$),
write $f(x)=j(x)^2 h(x)$ in $\F_q[x]$ with $h(x)$ squarefree,
so by the previous paragraph $\deg h > 0$.\footnote{
	At this point we know that $y^2=f(x)$
	has a multiplicity one component
	which
is geometrically integral, so it follows automatically as in
	Proposition~3.1 in~\cite{vg1} or Lemma~15 in~\cite{poonenstoll}
	that $\C_t$ has a $k$-rational divisor of degree~1.
	In fact this would suffice for our application,
	except in Section~\ref{BMOforX}, where it will be convenient
	to have the stronger result $\C_t(k) \not= \emptyset$.}
Any $\F_q$-point on $z^2=h(x)$ with $x$ finite and $j(x) \not=0$
gives rise to a nonsingular $\F_q$-point
on the affine curve $y^2=f(x)$,
which can be lifted by Hensel's Lemma to prove $\C_t(k)\not=\emptyset$.
Also, if there is a point at $\infty$ on
(the projective nonsingular model of) $z^2=h(x)$ over $\F_q$,
then there will be points on $\C_t$ at $\infty$.
The number of points on $z^2=h(x)$ with $j(x)=0$
is at most $2 \deg j = 6 - \deg h = 6 - (2g+2) = 4-2g$,
where $g$ is the geometric genus of $z^2=h(x)$.
Therefore by the Weil bound,
we automatically find the desired point if
	$$(q+1-2g\sqrt{q}) - (4-2g) > 0$$
for $g=0,1,2$.
These hold for $q \ge 17$, so it remains to
consider the cases $q=3,5,7,9,11,13$.
If $q \equiv 1 \pmod{4}$, then $\C_t$ has $k$-points at $\infty$.
For each remaining $q$,
we check by hand that for each possible residue class of $t$
(and for the case $T \equiv 0$ in~(\ref{model2})),
there is still a nonsingular affine point on $y^2=f(x)$.
\end{proof}

\subsection{Components of the Picard scheme}
\label{components}

Since $\C$ is a proper smooth surface over $\Q$,
it is projective over $\Q$, and hence also projective over $\PP^1$.
Moreover the geometric fibers of $\C \rightarrow \PP^1$ are integral.
Hence we may apply 9.3/1 from~\cite{blr}
to deduce that
\begin{enumerate}
	\item The relative Picard functor $\PIC_{\C/\PP^1}$
		is represented by a separated and locally smooth
		$\PP^1$-scheme.
	\item $\PIC_{\C/\PP^1} = \coprod_{n \in \Z} \PIC_{\C/\PP^1}^n$,
		where $\PIC_{\C/\PP^1}^n$ denotes the open and closed
		subscheme of $\PIC_{\C/\PP^1}$ consisting of all line bundles
		of degree~$n$.
	\item For each $n \in \Z$,
		$\PIC_{\C/\PP^1}^n$ is a torsor for $\PIC_{\C/\PP^1}^0$
		and is quasi-projective.
\end{enumerate}
Let $\Atilde:=\PIC_{\C/\PP^1}^0$ and let $\Xtilde:=\PIC_{\C/\PP^1}^1$.
For any $\PP^1$-scheme $\pi: Z \rightarrow \PP^1$,
let $Z_U$ denote $\pi^{-1}(U)$.
Then $\Atilde_U$ is an abelian $U$-scheme, by 9.4/4 from~\cite{blr}.
Resolution of singularities gives us $\PP^1$-schemes
$\pi_\A: \A \rightarrow \PP^1$ and $\pi_\X: \X \rightarrow \PP^1$ such that
\begin{enumerate}
	\item $\A$ and $\X$ are smooth and projective over $\Q$.
	\item $\A_U \isom \Atilde_U$ and $\X_U \isom \Xtilde_U$
		as $U$-schemes.
\end{enumerate}

We now verify that $\A$ and $\X$
satisfy the conditions of Theorem~\ref{main}.
Conditions~(1) and~(2) are clear from the construction.
The following lemma applied to the smooth fibers $\C_t$ shows that
$\pi_\X(\X(k_v)) \supset U(k_v)$,
and condition~(3) then follows from compactness of $\pi_\X(\X(k_v))$.

\begin{lemma}
\label{lichweich}
Let $k$ be a local field (i.e., a finite extension of $\R$, $\Q_p$, or
$\F_p((t))$), and let $\kbar$ denote a separable closure of $k$.
Let $X$ be a smooth projective geometrically integral
curve over $k$ of genus $g$,
and let $\overline{X} = X \times_k \kbar$.
Then there exists an element of degree~$g-1$ in $\Pic(\overline{X})$ which is
stable under $\Gal(\kbar/k)$.
\end{lemma}

\begin{proof}
This is a consequence of Tate local duality.
See Theorem~7 on p.~133 of~\cite{lichtenbaum}
for the case where $k$ is a finite extension of $\Q_p$,
or Section~4 of~\cite{poonenstoll} for the general case.
The real case is related to W.-D.~Geyer's
modern version of work of Weichold and Witt:
see~\cite{scheiderer}, Thm.~(20.1.5.1),  p.~221.
\end{proof}

We now prove condition~(4).
Suppose condition~(4) fails, so that
$\X$ has a closed point of odd degree.
Its image under $\pi$ is a a closed point of $\PP^1$
of odd degree.
Since $U$ contains all points of $\PP^1$ of odd degree,
this closed point corresponds to some $t \in U(k)$ with $[k:\Q]$ odd,
such that the fiber $\X_t$ has a closed point of
odd degree over $k$.
Since  $\C_t$ admits a 2-to-1 map to $\PP^1_k$, it possesses
a point in a quadratic extension of $k$, hence $\X_t=\PIC_{\C_t/k}^1$
possesses also such a point. Thus the principal homogeneous space
$\X_t=\PIC_{\C_t/k}^1$ of $\A_t=\PIC_{\C_t/k}^0$ is annihilated
by coprime integers, hence $\X_t$ is trivial.

On the other hand, we will show that $\X_t$ cannot be trivial
because of the following,
which is Corollary~12 in~\cite{poonenstoll} (the assertion
$c \in \Sha(A)$ is Lemma~\ref{lichweich} above).

\begin{lemma}
\label{cc}
Let $k$ be a global field (i.e., a finite extension of $\Q$ or $\F_p(t)$).
Let $C$ be a smooth projective geometrically integral curve over $k$.
Let $A$ be the canonically polarized jacobian of $C$,
and let $c$ be the class of $\PIC_{C/k}^{g-1}$ in $\HH^1(k,A)$.
Then $c \in \Sha(A)$, and the Cassels-Tate pairing satisfies
$\langle c,c \rangle = N/2 \in \Q/\Z$, where $N$ is the number of
places $v$
of $k$ for which $C$ has no $k_v$-rational divisor of degree $g-1$.
\end{lemma}

Let us apply this to the fiber $C=\C_t$.
Then $g=2$, $c$ is the class of $\X_t$,
and $N$ equals the number of real places of $k$,
by Lemma~\ref{localdivisors}.
Since $[k:\Q]$ is odd, $N$ is odd.
Thus $\langle c,c \rangle \not=0$, so $c \not=0$;
i.e., $\X_t$ is non-trivial in $\Sha(\A_t)$.
This contradiction proves condition~(4).

We next prove condition~(5).
The curve $\C_U$ is projective and smooth over $U$.
There is an associated $\Q$-morphism $f$ from
$U$ to the coarse moduli space $M_2$ of smooth
curves of genus~2, and the jacobian functor yields a
$\Q$-morphism $j : M_2 \rightarrow \A_{2}$
to the coarse moduli space of
principally polarized abelian surfaces over ${\bf Q}$.
By Torelli's theorem,
the map $j$ is injective over algebraically
closed fields~\cite[p.~143]{GIT}.
Thus to prove~(5), it is enough to
show that $f$ is not constant.
Suppose $f$ is constant.
Let $R$ be the completion of the local ring of $\PP^1$
at the closed point where
$\C/\PP^1$ has bad reduction,
and let $K$ be the fraction field of $R$.
Since $f$ is constant,
there exists a finite field extension $L/K$ such that
$\C \times_{\PP^1} L$ is $L$-isomorphic to a constant curve.
The latter has good reduction over the ring of integers of $L$,
but for a stable curve the property of having good
reduction does not depend on the field extension,
so the stable curve $\C \times_{\PP^1} R$ must be smooth.
This is a contradiction.

Finally we must show that the generic fiber of $\A \rightarrow \PP^1$
is absolutely simple.
Let $\Ctilde$ denote the (good) mod~3 reduction
of the fiber $\C_{-1}$, and let $J$ be the jacobian of $\Ctilde$.
Calculating $\#\Ctilde(\F_3)=1$ and $\#\Ctilde(\F_9)=15$,
we find that the characteristic polynomial of Frobenius
for $J$ is $x^4 -3 x^3 +7 x^2 -9 x +9$.
A root $\alpha$ generates a quartic field $L$ having a unique
non-trivial proper subfield $F$ (of degree~2).
The $L/F$-conjugate of $\alpha$ is $3\alpha^{-1}$,
and $3\alpha^{-1}/\alpha$ is not a root of unity,
so $\alpha^n \not\in F$ and $\alpha^n$ is of degree~4 for all $n \ge 1$.
Hence $J$ is absolutely simple.
To deduce from this that the generic fiber of $\A \rightarrow \PP^1$
is absolutely simple,
it suffices to apply the following well-known result twice:

\begin{lemma}
Let $\A$ be an abelian scheme over a discrete valuation ring $R$.
If the special fiber is absolutely simple,
then so is the generic fiber.
\end{lemma}

\begin{proof}
If the generic fiber is not absolutely simple,
then there exists a finite extension $K'$ of the fraction field $K$ of $R$
and an isogeny $\A \times_R K' \rightarrow B \times_{K'} C$,
where $B$ and $C$ are non-trivial abelian varieties over $K'$.
Let $R'$ be a discrete valuation ring in $L$ containing $R$,
and let $k'/k$ be the corresponding extension of residue fields.
By~\cite{serretate}, the N\'eron models $\B$ and $\C$
of $B$ and $C$ are abelian schemes over $R'$.
The universal property of the N\'eron model
extends the isogeny to the N\'eron models:
$\A \times_R R' \rightarrow \B \times_{R'} \C$.
Taking special fibers shows that $\A \times_R k'$
is not simple, so $\A \times_R k$ is not absolutely simple.
\end{proof}

This completes the proof of Theorem~\ref{main}.

\subsection{Review of the Brauer-Manin obstruction}
\label{review}

Manin~\cite{maninICM} in 1970
introduced what is now called the Brauer-Manin obstruction.
Now suppose $X$ is a smooth proper irreducible variety
over a number field $k$.
Let $\Ak$ denote the ad\`ele ring of $k$.
If $v$ is a place of $k$,
then we may define a ``local evaluation pairing''
	$$\ev_v: \Br(X) \times X(k_v) \rightarrow \Q/\Z,$$
by letting $\ev_v(\alpha,x) \in \Q/\Z$ be the invariant
of $i_x^\ast \alpha \in \Br(k_v)$
where $i_x: \Spec k_v \rightarrow X$
corresponds to $x \in X(k_v)$.
Given $\alpha$,
the properness of $X$ implies
that the map $x \mapsto \ev_v(\alpha,x)$ is zero for almost all $v$,
so by summing over all $v$ we obtain a global pairing
	$$\ev: \Br(X) \times X(\Ak) \rightarrow \Q/\Z,$$
continuous in the second argument,
such that $\ev(\alpha,x)=0$ if $\alpha$ comes from $\Br(k)$
or if $x \in X(k)$.
One says that there is a Brauer-Manin obstruction for $X$
if for every $x \in X(\Ak)$ there exists $\alpha \in \Br(X)$
with $\ev(\alpha,x) \not= 0$.

It is conjectured that for any smooth, proper, (geometrically)
rational surface $X$
over a number field $k$, the Brauer-Manin obstruction
is the only obstruction to the Hasse principle;
i.e., that if $X(k)=\emptyset$, it is because
$\{x \in X(\Ak): \ev(\alpha,x)=0 \text{ for all } \alpha \in \Br(X)\}
	= \emptyset$.
This appeared as a question in~\cite{cts},
and as a conjecture in~\cite{ctks}, once many examples were available.
Any rational surface over $k$ is $k$-birational either
to a conic bundle over a conic, or to a Del Pezzo surface
of degree~$d$, $1 \leq d \leq 9$. Del Pezzo surfaces of degree
at least~5 satisfy the Hasse principle. The Del Pezzo surfaces
of degree~4 are the smooth intersections of two quadrics in $\PP^4$,
and those of degree~3 are the smooth cubic surfaces in $\PP^3$.
There exists theoretical evidence for the conjecture in the case of conic
bundles and Del Pezzo surfaces of degree~4,
and numerical evidence for diagonal cubic surfaces.

\subsection{The evaluation pairing and the Cassels-Tate pairing}
\label{casselstate}

It is possible to define the Cassels-Tate pairing in terms
of the pairing $\ev$ defined above.
Suppose $X$ and $Y$ are torsors of an abelian variety $A$ over $k$,
representing elements of $\Sha(A)$.
We may identify $A$ with $\PIC^0_{X/k}$,
and then $Y$ corresponds to an element of $\HH^1(k,\Pic^0 {\overline X})$,
which may be mapped to an element $\beta \in \HH^1(k,\Pic {\overline X})$.
The Leray spectral sequence~\cite[III.1.18(a)]{milneetale}
for $X \rightarrow \Spec k$ gives rise to an exact sequence
	$$\ker( \Br(X) \rightarrow \Br({\overline X}))
	\longrightarrow \HH^1(k,\Pic {\overline X})
	\longrightarrow \HH^3(k,\G_m) = 0,$$
the last equality holding because $k$ is a number field.
We pick $x \in X(\Ak)$, and pick $\alpha \in \Br(X)$ mapping to $\beta$
in the sequence above.
Then the Cassels-Tate pairing $\langle X,Y \rangle$ equals
$\ev(\alpha,x)$, which does not depend on the choices made.
(See Theorem~41.24 in~\cite{manincubic} for the case of genus~1 curves,
or the ``homogeneous space definition'' of the pairing in~\cite{poonenstoll}.)

\subsection{The Brauer-Manin obstruction for $\X$}
\label{BMOforX}

It is conjectured in~\cite[Conj.~2]{ct}
that the Brauer-Manin obstruction
to the existence of zero-cycles of degree~1
(the definition is an obvious extension of the one given for rational points)
on smooth projective irreducible varieties over number fields
is the only one.
In this section we show that the non-existence of zero-cycles of degree~1
on our variety $\X$ can
indeed be explained by this obstruction; moreover it can be explained
by the obstruction attached to a
single element of $\Br(\X)$.

\begin{lemma}
\label{H3kC}
Let $k$ be a number field,
and let $C$ be a smooth projective geometrically integral curve over $k$.
Then $\HH^3(k(C),\G_m)=0$.
\end{lemma}

\begin{proof}
If $F$ is any field of characteristic zero,
$\mu_\infty$ denotes the group of roots of unity in $\overline{F}$,
and $Q:={\overline{F}}^\ast/\mu_\infty$, then $Q$ is uniquely divisible,
so $\HH^i(F,Q)=0$ for $i>0$,
and we deduce that $\HH^3(F,\G_m)=\HH^3(F,\mu_\infty)=\HH^3(F,\Q/\Z(1))$.

For real places $v$ of the number field $k$,
let $k_v':=k_v \cap \kbar \subset \kbar_v$ denote
the corresponding real closure of $k$.
Theorem~$1'$ of~\cite{jannsen} states that
$\HH^3(k(C),\Q/\Z(r)) \isom \bigoplus_{v \real} \HH^3(k_v'(C),\Q/\Z(r))$
for $r \not= 2$.
Taking $r=1$, we find that
$\HH^3(k(C),\G_m) \isom \bigoplus_{v \real} \HH^3(k_v'(C),\G_m)$.

For $v$ real,
let $G_v'=\Gal(\kbar(C)/k_v'(C)) \isom \Gal(\kbar/k_v') \isom \Z/2$.
The spectral sequence
	$$E_2^{p,q} := \HH^p(G_v',\HH^q(\kbar(C),\G_m))
		\implies \HH^{p+q}(k_v'(C),\G_m)$$
together with the fact (Th\'eor\`eme~1.1 in~\cite{grothbrauer})
that $\HH^q(\kbar(C),\G_m)=0$ for $q>0$
implies that
	$$\HH^3(k_v'(C),\G_m) = \HH^3(G_v',\kbar(C)^\ast).$$
By periodicity, the latter is $\HH^1(G_v',\kbar(C)^\ast)$,
which is zero by Hilbert's Theorem 90.
\end{proof}

\begin{cor}
\label{directlimit}
Let $k$ be a number field,
and let $C$ be a smooth projective geometrically integral curve over $k$.
Then $\varinjlim \HH^3_\et(V,\G_m)=0$,
where the direct limit is over the dense open subschemes $V$ of $C$,
with respect to the restriction maps.
\end{cor}

\begin{proof}
The direct limit equals $\HH^3(k(C),\G_m)$, by
a result of Grothendieck.
See Corollaire~5.9 in SGA~4~VII (Expos\'e VII : Site et Topos \'etale d'un
sch\'ema, par A. Grothendieck)~\cite{sga4},
or Lemma~1.16 and Remark~1.17(a) in~\cite{milneetale}.
\end{proof}

\begin{rem}
With more work one can show that if $U \subset C$ is an affine open subset,
then $\HH^3_\et(U,\G_m)$ is finite, but not necessarily zero.
Moreover, there exists a dense open subset $U \subset C$
such that $\HH^3_\et(U',\G_m)=0$ for all open subsets $U' \subseteq U$.
On the other hand, $\HH^3_\et(C,\G_m)$ is always infinite.
\end{rem}

We now return to our situation,
with $\C$ and $U$ as in Section~\ref{relativecurve}
and $\X$ as in Section~\ref{components}.
Let $V \subseteq U$ be a dense open subset.
Recall that $\X_V=\PIC^1_{\C_V/V}$.
Let $\iota: \C_V \rightarrow \X_V$ be the natural inclusion of $V$-schemes.
We obtain a commutative diagram of exact sequences
of \'etale sheaves over $V$
$$\begin{CD}
	&& \PIC^0_{\X_V/V} @>>> \PIC_{\X_V/V}	\\
	&& @V{\iota^\ast}VV @V{\iota^\ast}VV 	\\
	0 @>>> \PIC^0_{\C_V/V} @>>> \PIC_{\C_V/V} @>>> \Z @>>> 0
\end{CD}$$
and an argument similar to
the one of
Proposition~6.9 on p.~118 of~\cite{GIT}
shows that the first vertical map is an isomorphism.
Taking cohomology we obtain
$$\begin{CD}
	&& \HH^1_\et(V,\PIC^0_{\X_V/V}) @>>> \HH^1_\et(V,\PIC_{\X_V/V})	\\
	&& @V{\iota^\ast}VV @V{\iota^\ast}VV 	\\
	\Z @>>> \HH^1_\et(V,\PIC^0_{\C_V/V}) @>>> \HH^1_\et(V,\PIC_{\C_V/V})
\end{CD}$$
The image of $1 \in \Z$ in $\HH^1_\et(V,\PIC^0_{\C_V/V})$
is the class $\gamma$ of the torsor $\X_V$ of $\PIC^0_{\C_V/V}$.
(The analogous result for the generic fibers, with $V$ replaced by
its generic point $\eta:=\Spec \Q(t)$,
follows from a cochain calculation in Galois cohomology.
The claim for $V$ follows, using functoriality and
applying the injectivity of
$\HH^1_\et(V,\B) \rightarrow \HH^1(\eta,\B_\eta)$
for abelian $V$-schemes $\B$ to the case $\B=\PIC^0_{\C_V/V}$.)

Let $\beta$ be the image of $(\iota^\ast)^{-1}(\gamma)$
under $\HH^1_\et(V,\PIC^0_{\X_V/V}) \rightarrow \HH^1_\et(V,\PIC_{\X_V/V})$.
Then $\beta$ maps to $0$ in $\HH^1_\et(V,\PIC_{\C_V/V})$.
The Leray spectral sequences for $\pi_{\X_V} :\X_V \rightarrow V$
and $\pi_{\C_V}: \C_V \rightarrow V$ yield a
commutative diagram of exact sequences
\begin{equation}
\label{doubleLeray}
\begin{CD}
	\Br V @>>>
	\ker \left[ \Br \X_V \rightarrow \HH^2_\et(V,\RR^2 \pi_{\X_V \ast}
\G_m) \right]
	@>>> \HH^1_\et(V,\PIC_{\X_V/V}) @>>> \HH^3_\et(V,\G_m)	\\
	@| @VVV @VVV	\\
	\Br V @>>>
	\ker \left[ \Br \C_V \rightarrow \HH^2_\et(V,\RR^2\pi_{\C_V \ast}
\G_m) \right]
	@>>> \HH^1_\et(V,\PIC_{\C_V/V})
\end{CD}
\end{equation}
By Corollary~\ref{directlimit},
after shrinking $V$ and restricting if necessary,
we may assume that $\beta \in \HH^1_\et(V,\PIC_{\X_V/V})$
maps to zero in $\HH^3_\et(V,\G_m)$,
so that $\beta$ lifts leftwards
to an element $\alpha \in \Br \X_V$.
Since $\beta$ maps to 0 in $\HH^1_\et(V,\PIC_{\C_V/V})$,
we may adjust $\alpha$ by an element in the image of $\Br V$ if necessary,
to assume that $\alpha$ restricts to 0 in $\Br \C_V$.

Suppose that $t \in V(\Q_p)$ for some finite or infinite prime $p$,
and consider the fiber $\X_t$,
which is a torsor of an abelian variety over $\Q_p$.
The restriction of $\gamma$ to an element $\gamma_t$
of $\HH^1(\Q_p,\PIC^0_{\X_t/\Q_p})$
represents the torsor $\PIC^1_{\C_t/\Q_p}$,
which is trivial by Lemma~\ref{lichweich}.
Then the corresponding restriction $\beta_t \in \HH^1(\Q_p,\PIC_{\X_t/\Q_p})$
is 0, so the analogue of the top row of~(\ref{doubleLeray})
for $\X_t \rightarrow \Spec \Q_p$
shows that the restriction $\alpha_t$ comes from $\Br \Q_p$.
In particular, $\ev_p(\alpha_t,x)$ is a constant function of
$x \in \X_t(\Q_p)$.
For any finite $p$,
Lemma~\ref{localdivisors} gives us a point $x \in \C_t(\Q_p)$,
but $\alpha$ restricts to $0 \in \Br \C_t$,
so the constant function must be identically zero.
Thus for all finite $p$,
$\ev_p(\alpha,x)=0$ for
all $x \in \X_V(\Q_p)$.
By~\cite[Th\'eor\`eme 2.1.1]{harari},
this implies that $\alpha$ extends (uniquely)
to an element $\tilde{\alpha} \in \Br \X$.

By continuity, $\ev_p(\tilde{\alpha},x)=0$ for all $x \in \X(\Q_p)$,
for any finite $p$.
The evaluation of $\tilde{\alpha}$
on an adelic point $x=(x_v) \in \X_V(\A_\Q)$
depends only on $x_\infty$.
The union of the points $x_\infty \in \X(\R)$
lying above {\em rational} points $t \in V(\Q)$
is dense in $\X(\R)$,
and by Section~\ref{casselstate},
$\ev(\tilde{\alpha},x)$ for any such $x$ equals
the value of the Cassels-Tate pairing $\langle \X_t,\X_t \rangle$,
which is $1/2$.
Hence by continuity $\ev(\tilde{\alpha},x)=1/2$ for all $x \in \X(\A_\Q)$,
and we deduce {\em a posteriori} that
$\ev_\infty(\tilde{\alpha},x)=1/2$ for all $x \in \X(\R)$.

The same argument which
showed that
$\ev_p(\tilde{\alpha},x)=0$ for all $x \in \X(\Q_p)$
shows that if $k_v$ is any nonarchimedean local field containing $\Q$
and $x \in \X(k_v)$, then
$\ev_v(\tilde{\alpha},x)=0$.

{}From the above two paragraphs,
we see that for any zero-cycle $z_v$ on $\X \times_\Q \Q_v$,
the element $\ev_v(\tilde{\alpha},z_v)$
equals $0 \in \Q/\Z$ if $v$ is a finite place,
and equals $d_v/2 \in \Q/\Z$ if $v$ is the real place,
where $d_v$ is the degree of the zero-cycle $z_v$ (over $\R$).
If the $z_v$ for all $v$ arise from a zero-cycle $z$ on $\X$ of odd degree,
then summing over $v$ shows that $\ev(\tilde{\alpha},z)=1/2$,
which is impossible.
By definition, this means that $\tilde{\alpha}$
gives a Brauer-Manin obstruction to the existence of such a zero-cycle $z$.

\section{A family of genus~1 curves}
\label{genus1family}

We now prove Theorem~\ref{genus1}.

\subsection{Cubic surfaces violating the Hasse principle}
\label{violation}

Swinnerton-Dyer~\cite{swinnertondyer}
disproved a conjecture of Mordell~\cite{mordell}
by constructing a smooth cubic surface $V$ in $\PP^3$ over $\Q$
``violating the Hasse principle'';
i.e., such that $V$ has points over each completion of $\Q$,
but not over $\Q$.
Soon afterwards,
Cassels and Guy~\cite{casselsguy} gave the {\em diagonal} cubic surface
	$$5 x^3 + 9 y^3 + 10 z^3 + 12 t^3 = 0$$
violating the Hasse principle.

Manin~\cite{maninICM} used the Brauer-Manin obstruction
to explain Swinnerton-Dyer's counterexample.
(See also~\cite{manincubic}.)
Much later~\cite{ctks} explained
the Cassels-Guy example from this point of view,
and gave a very explicit algorithm for computing
the Brauer-Manin obstruction for diagonal cubic surfaces.

We fix a number field $k$
and a smooth cubic surface $V \subset \PP^3_k$
for which the Hasse principle fails.

\subsection{Lefschetz pencils}

By~2.5.2 in~\cite{katz}, the embedding $V \hookrightarrow \PP^3_k$
is a ``Lefschetz embedding.''
This implies that the statements in the rest of this paragraph hold
for a sufficiently generic choice of a $k$-rational line $L$
in the dual projective space $\PPdual^3_k$.
(``Sufficiently generic'' here means for $L$ corresponding to
points outside a certain Zariski closed subset of the Grassmannian.)
Let $L'$ denote the
{\em axis} of $L$, that is, the line in $\PP^3_k$
obtained by intersecting two hyperplanes
in the family given by $L$.
Blowing up $V$ at the scheme-theoretic intersection $V \cap L'$
results in a smooth projective variety $V'$ over $k$
isomorphic to the reduced
subvariety of $V \times L$ whose geometric points
are the pairs $(v,H)$ where $v \in V$ is on the hyperplane $H \subset \PP^3_k$
corresponding to a point of $L$.
Each fiber of $V' \rightarrow L \isom \PP^1_k$
is a proper geometrically integral
curve of arithmetic genus~1, and if it is singular, there is only
one singularity and it is a node.
The generic fiber is smooth.

\subsection{Local points in the pencil}
\label{localpoints}

Let $k_v$ denote the completion of $k$ at a non-trivial place $v$.
Let $\OO_v$ be the ring of integers, let $\mm_v$ be its maximal ideal,
and let $\F_v=\OO_v/\mm_v$.

Variants of the following lemma
have appeared in various places in the literature (e.g. \cite{ctssd},
\cite[2.1]{ct}). The key point is that {\it all} the
geometric fibers of $V' \rightarrow \PP^1$ are
integral.

\begin{lemma}
\label{surjectivity}
There exists a finite set $S$ of places of $k$
such that for any  $v \not \in S$ and any finite extension $K$ of $k_v$
the map $V'(K) \rightarrow \PP^1(K)$ is surjective.
\end{lemma}

\begin{proof}
The morphism $f: V' \rightarrow \PP^1_k$ is proper and flat.
Combining Theorem~11.1.1 and Theorem~12.2.4(viii) of~\cite{EGA4},
we see that there exists a
a finite set $S$ of places of $k$, with associated ring of
$S$-integers $\OO$,
such that $f$ extends to a proper flat morphism
$f' : \V' \rightarrow \PP^1_\OO$,
all fibers of which
are geometrically integral.
We may assume that $S$ contains all the real places
and none of the complex places.
The desired surjectivity is automatic at complex places,
so it remains to prove the surjectivity for finite $v \not\in S$.

Let $v$ be a place not in $S$, let $K/k_v$ be a finite extension
with associated valuation $w$,
let $\OO_w$ be the ring of integers of $K$ and let
 $\F_w$ be the residue field.
Let $q=\#\F_w$.
If $t \in \PP^1(\OO_w)=\PP^1(K)$
then the special fiber of the fiber $\V'_t$
is a geometrically integral, proper curve of arithmetic genus~1.
The number of smooth $\F_w$-points on this special fiber
is at least $(\sqrt{q}-1)^2 > 0$ (Hasse)
if the special fiber is smooth,
and is equal to $q-1$, $q$, or $q+1$ if not;
in any case there is at least one.
Hensel's Lemma then shows that the generic fiber of $\V'_t$
has a $K$-point, as desired.
\end{proof}

Let $U'$ be the largest open subscheme of $\PP^1_k$
over which $V'$ is smooth.

\begin{lemma}
\label{opensubset}
For each completion $k_v$ of $k$,
the image of $V'(k_v) \rightarrow \PP^1(k_v)$
contains a nonempty open subset $W_v$ (in the $v$-adic topology).
\end{lemma}

\begin{proof}
Recall that $V(k_v) \not= \emptyset$ by choice of $V$.
Since $V'$ and $V$ are birational smooth projective varieties,
$V'(k_v) \not= \emptyset$ too,
and in fact $V'(k_v)$ is Zariski dense in $V'$.
In particular we can find $P \in V'(k_v)$ mapping into $U'$.
The image of $V'(k_v) \rightarrow \PP^1(k_v)$ will then
contain a neighborhood of the image of $P$.
\end{proof}

We may choose $S$ as in Lemma~\ref{surjectivity}
to contain no complex places.
For $v \in S$ choose $W_v$ as in Lemma~\ref{opensubset}.
We may assume that there exists $v \in S$
for which $W_v \subset U'(k_v)$.
Then the fibers $V'_t$ for $t \in W_v$ are smooth.

\subsection{Base change of the pencil}

\begin{lemma}
\label{adelicimage}
Let $S$ be a finite set of non-complex places of $k$,
let $T$ be a finite set of closed points of $\PP^1_k$,
and for each $v \in S$, let $W_v$ be a nonempty open subset
of $\PP^1(k_v)$.
Then there exists a non-constant $k$-morphism
$f: \PP^1_k \rightarrow \PP^1_k$
such that $f$ is smooth (i.e., unramified) above the points of $T$,
and $f(\PP^1(k_v)) \subseteq W_v$ for all $v \in S$.
\end{lemma}

\begin{proof}
First we show that for each $v \in S$,
we can find a non-constant $f_v$ defined over $k_v$
such that $f_v(\PP^1(k_v)) \subseteq W_v$.
Choose an affine point $\lambda \in W_v$.
If $v$ is real, we may take $f_v(x)=\lambda+(x^2+n)^{-1}$ for $n \gg 0$.
If $v$ is finite, then we may choose $g \in k_v(x)$
such that $g$ maps $\{0,1,\infty\}$ to $\lambda$,
and then let $f_v(x)=g(x^n)$ for $n=q^k(q-1)$
where $q=\#\F_v$ and $k \gg 0$.

Preceding each $f_v$ with an arbitrary non-constant rational function
of appropriate degree, we may assume that $\deg f_v$ is the same
for all $v \in S$.
Let $B_d$ denote the open affine subset
of $(a_0,\dots,a_d,b_0,\dots,b_d) \in \Aff^{2d+2}$
for which the homogeneous polynomials
$\sum a_i X^i Y^{d-i}$ and $\sum b_i X^i Y^{d-i}$
have no non-trivial common factor.
We have a morphism $B_d \times \PP^1 \rightarrow \PP^1$
which constructs the rational function of degree~$d$
which is the quotient of the two polynomials,
and then evaluates it at the point in $\PP^1$.
Choose $b_v \in B_d(k_v)$ representing $f_v$.
By compactness of $\PP^1(k_v)$ in the $v$-adic topology,
any point of $B_d(k_v)$ sufficiently close to $b_v$
represents a rational function still mapping $\PP^1(k_v)$ into $W_v$.

Weak approximation gives us a point $b \in B_d(k)$
close enough $v$-adically to $b_v$ for each $v \in S$
so that the corresponding rational function $f$ over $k$
maps $\PP^1(k_v)$ into $W_v$ for $v \in S$.
Moreover $f$ will be smooth above the points of $T$
provided that $b$ avoids a certain closed subset of $B_d$
of positive codimension, so this is easily arranged.
\end{proof}

Apply Lemma~\ref{adelicimage} to obtain $f$
for the $S$ and the $W_v$ at the end of Section~\ref{localpoints},
and with $T=\PP^1_k \setminus U'$.
We let $\pi_\X: \X \rightarrow \PP^1_k$ be the base extension
of $V' \rightarrow \PP^1_k$ by $f$, so that the following is
a cartesian square:
$$\begin{CD}
	\X		@>>> V'		\\
	@V{\pi_\X}VV	@VVV	\\
	\PP^1_k		@>f>> \PP^1_k
\end{CD}$$
Each factor is smooth over $k$,
and above each point of the lower right $\PP^1_k$
at least one of the two factors is smooth
by  choice of $T$, so $\X$ is smooth over $k$.
The generic fiber of $V' \rightarrow \PP^1_k$ is geometrically integral,
so $\X$ is geometrically integral.
Since $V'$ is projective over $k$, so is $\X$.

Let $U=f^{-1}(U')$.
By choice of $f$, $f(\PP^1(k)) \subset U'$.
Hence $U$ is an open subscheme of $\PP^1_k$ containing $\PP^1(k)$,
and $\X_U \rightarrow U$ is a proper smooth family
of geometrically integral curves of genus~1.
We may construct $\E$ as the minimal proper regular model of
the jacobian of the generic fiber of $\X \rightarrow \PP^1_k$.
Then $\E_U \isom \PIC_{\X_U/U}^0$ is an abelian $U$-scheme.

If $v$ is a place of $k$, and $t \in \PP^1(k_v)$,
then $\X_t \isom V'_{f(t)}$ has a $k_v$-point by choice of $f$.
The property of having a rational point is
a birational invariant of smooth projective varieties,
and $V(k)=\emptyset$, so $V'(k)=\emptyset$.
Since $\X$ maps to $V'$, $\X(k)=\emptyset$ too.

It remains to show that the $j$-invariant of
$\pi_\X: \X \rightarrow \PP^1_k$
is non-constant,
or equivalently that the $j$-invariant of $\pi: V' \rightarrow \PP^1_k$
is non-constant.
Euler-Poincar\'e characteristic
computations~\cite[Prop.~(11.4), p.~97]{bpv}
show that the fibration $V' \rightarrow \PP^1$
has exactly 12 nodal geometric fibers.
Thus the $j$-invariant of the generic fiber has poles on $\PP^1$,
so it cannot be constant.

\subsection{Generic Cassels-Tate pairing}

Let $k$ be a number field,
let $U \subset \PP^1_k$ be a dense open subscheme of $\PP^1_k$,
and let $\A$ be an abelian $U$-scheme.
Assume that $\A \rightarrow U$ is projective,
so that the dual abelian scheme $\A' \rightarrow U$ exists~\cite{grothpicard}.
Now let $\X_U$ be an $\A$-torsor over $U$
such that each fiber $\X_t$ for $t \in U(k)$ represents
a nonzero element of $\Sha(\A_t)$.
Is it then automatic that there exists a $\A'$-torsor $\Y_U$ over $U$
(or at least over some dense open subscheme)
such that for all $t \in U(k)$,
(the class of) $\Y_t$ is in $\Sha(\A'_t)$
and the Cassels-Tate pairing $\langle \X_t,\Y_t \rangle$
gives a nonzero value in $\Q/\Z$ which does not depend on $t$?
The special case where the families $\A$ and $\X$ are split
is the well-known conjecture
that the Cassels-Tate pairing is nondegenerate.

For our first example, the $\X_U$ in Section~\ref{abeliansurfaces}
built out of the $\PIC^1$ of the relative curve of genus~2,
we can identify $\A$ with $\A'$ and take $\Y_U=\X_U$,
because we showed $\langle \X_t,\X_t \rangle = 1/2 \in \Q/\Z$
for all $t \in U(\Q)$.

For the second example,
the one from the Lefschetz pencil in a cubic surface $V$ over $k$,
we can construct $\Y_U$, if we assume (as is conjectured)
that the failure of the Hasse principle for $V$
is due to a Brauer-Manin obstruction, as we now explain.

The cokernel of
$\Br(k) \rightarrow \Br(V)$ for a cubic surface $V$ is finite,
but need not be cyclic:
see~\cite{swinnertondyer2} for the possibilities.
Nevertheless one has the following:

\begin{lemma}
\label{oneelement}
Let $V$ be a smooth cubic surface in $\PP^3$ over a number field $k$.
Assume that the Hasse principle for $V$ fails,
and fails because of a Brauer-Manin obstruction for $V$.
Then there exists $\alpha \in \Br(V)$
such that $\ev(\alpha,x)$ is a nonzero constant
independent of $x \in V(\Ak)$.
\end{lemma}

\begin{proof}
Let $G$ be the cokernel of $\Br(k) \rightarrow \Br(V)$,
and let $\Ghat = \Hom(G,\Q/\Z)$.
These are finite, and are given the discrete topology.
Let $S$ be the image of the continuous map $\phi: V(\Ak) \rightarrow \Ghat$
induced by $\ev$.

We claim that
\begin{equation}
\label{coset}
	x,y,z \in S \implies y+z-x \in S.
\end{equation}
Because $\phi$ is continuous and $V$ is smooth,
given $x,y,z \in S$,
we may choose $P=(P_v) \in V(\Ak)$ with $\phi(P)=x$
and similarly $Q$ and $R$ giving $y$ and $z$,
so that
for each $v$, $P_v$, $Q_v$ and $R_v$ are not collinear,
and the plane through them intersects $V$
in a nonsingular genus 1 curve $C_v$.
Applying Riemann-Roch to $C_v$ yields a point $T_v \in C(k_v)$
linearly equivalent to the divisor $Q_v+R_v-P_v$ on $C_v$.
As in~\cite{lichtenbaum},
the pairing $\Br(C_v) \times C_v(k_v) \rightarrow \Br(k_v) \subseteq \Q/\Z$
extends to a pairing
$\Br(C_v) \times \Div(C_v) \rightarrow \Q/\Z$
which induces a pairing
$\Br(C_v) \times \Pic(C_v) \rightarrow \Q/\Z$.
For $\beta \in \Br(V)$,
$\ev(\beta,T)$ can be obtained by restricting $\beta$ to $C_v$,
pairing with $T_v$ and summing over $v$.
It follows that $\phi(T)=\phi(Q)+\phi(R)-\phi(P)$,
so $y+z-x \in S$.

It follows from~(\ref{coset}) and $S \not= \emptyset$
that $S$ is a coset of a subgroup $H$ of $\Ghat$.
But $S$ cannot be a subgroup, because $0 \in S$ would contradict
the existence of a Brauer-Manin obstruction.
Hence we may pick $g \in G \isom \Hom(\Ghat,\Q/\Z)$
annihilating $H$ but not $S$.
For any lift $\alpha \in \Br(V)$ of $g$,
$\ev(\alpha,x)$ is a nonzero constant for $x \in V(\Ak)$.
\end{proof}

\begin{rems} $\nichts$
\begin{enumerate}
\item Let $Y$ denote the Grassmannian of lines in $\PP^3$ over $k$.
Using correspondences and using the triviality of $\Br(Y)/\Br(k)$,
one can prove
	$$x,y \in S \implies -x-y \in S,$$
which improves~(\ref{coset}).
Thus $\alpha$ may be taken so that its image in $\Br(V)/\Br(k)$ has order~3.
(This follows also from the calculation of possibilities for $\Br(V)/\Br(k)$
in~\cite{swinnertondyer2},
together with Corollary~1 of~\cite{swinnertondyer2},
which forces $\Br(V)/\Br(k)$ to have odd order.)
Class field theory shows that $3 \Br(k) = \Br(k)$,
and it follows formally that the order~3 element of $\Br(V)/\Br(k)$
can be lifted to an order~3 element of $\Br(V)$.

\item
One can prove a similar result, namely
$x,y,z \in S \implies -x-y-z \in S$,
for del Pezzo surfaces $V$ of degree~4
(smooth intersections of two quadrics in $\PP^4$).
In this case $\Br(V)/\Br(k)$ is killed by~2.
(See p.~178 in~\cite{manincubic} or the proof of Proposition~3.18
in~\cite{ctssd}.)
If there is a Brauer-Manin obstruction, then it can be explained
by a single element $\alpha \in \Br(V)$ of order~2.
(If $k$ has real places, $2 \Br(k) \neq \Br(k)$,
so the argument in the previous remark
needs to be modified in order to prove this;
one must make use of the assumption that $V(k_v) \neq \emptyset$
for all real places $v$.)

\item
Suppose that $V$ is a del Pezzo surface of degree~4 over $k$ as above,
with a Brauer-Manin obstruction given by $\alpha \in \Br(V)$
of order~2.
A result of Amer~\cite{amer} and Brumer~\cite{brumer}
then implies that
$V(L)=\emptyset$ for any finite extension $L$ of odd degree over $k$.
According to the conjecture in Section~\ref{review},
there should be a Brauer-Manin obstruction for $V \times_k L$.
Does the same $\alpha$ yield an obstruction over $L$?

\item
If one replaces $V$ by a smooth cubic surface in the previous question,
and the Amer-Brumer result by the conjecture that the existence
of a point of degree prime to~3 on $V$ implies the existence of
a $k$-point,
then one is led to ask the analogous question for $V$,
for extensions $L$ with $\gcd([L:k],3)=1$.

\item There is a geometric application of the Amer-Brumer result
that leads to a question about Brauer-Manin obstructions for
Weil restriction of scalars.
Let $L$ be a finite separable field extension of $k$
and $W$ be a quasi-projective variety over $L$.
If $R:=\Res_{L/k}W$ is the Weil restriction of scalars,
then there is a natural $L$-morphism $R_L:=R \times_k L \rightarrow W$.
Now let $V$ be a del Pezzo surface of degree~4 over $k$,
and apply the above to $W:=V_L=V \times_k L$.
We obtain a morphism $R_L \rightarrow V_L$.
If $[L:k]$ is odd,
then application of the Amer-Brumer result for $V$
to the extension of function fields $L(R)/k(R)$,
yields a $k$-rational map $R \dashrightarrow V$.
A Brauer-Manin obstruction to the existence of a $k$-point on $V$
can be pulled back to obtain a Brauer-Manin obstruction for $R$ over $k$.
Hence we are led to the following question, to which a positive answer
would imply a positive answer to the question in remark~3 above:

Let $L/K$ be a finite extension of number fields.
Let $V$ be a smoooth projective geometrically integral variety over $L$,
and let $R=\Res_{L/k}(V)$.
It is clear that $V$ has points over all completions of $L$
if and only if $R$ has points over all completions of $k$.
Suppose that this is the case.
Then is it true that
there exists a Brauer-Manin obstruction for $V/L$ if and only if
there exists a Brauer-Manin obstruction for $R/k$?

If it is true that the Brauer-Manin obstruction is the only one
for geometrically rational varieties,
then the answer must be yes for such varieties.
For arbitrary $V$, one direction can be proved without too much work:
if there is a Brauer-Manin obstruction for $V/L$,
then because of the natural map $R_L \rightarrow V$ described above,
there is a Brauer-Manin obstruction for $R_L/L$,
and compatibility of the corestriction map
$\Br(R_L) \rightarrow \Br(R)$ with pullback
shows that there is a Brauer-Manin obstruction for $R/k$.
\end{enumerate}
\end{rems}

We return to the notation of earlier subsections
of Section~\ref{genus1family}.
In particular, $V$ is a smooth cubic surface violating the Hasse principle.
By our assumption (the conjecture), there is a Brauer-Manin obstruction.
Recall that we are trying to show
that the lack of rational points on $V$ can be explained
by the Cassels-Tate pairing.

Choose $\alpha$ as in Lemma~\ref{oneelement}.
Let $g$ be the composite morphism
$\X \rightarrow V' \rightarrow V$,
and define $\alpha' := g^\ast \alpha \in \Br(\X)$,
which gives a Brauer-Manin obstruction for $\X$.
We take $U$ to be the largest open subscheme of $\PP^1_k$
over which $\pi_\X: \X \rightarrow \PP^1_k$ is smooth,
and let $\X_U := \pi_X^{-1}(U)$ as usual.
We have an exact sequence of \'etale sheaves of groups
	$$0 \rightarrow \PIC^0_{\X_U/U} \rightarrow \PIC_{\X_U/U}
		\rightarrow \Z \rightarrow 0,$$
and $\HH^1_\et(U,\Z)=0$ by Proposition~3.6(ii) in~\cite{artin}.
This explains the bottom row of
\begin{equation}
\label{bigLeray}
\begin{CD}
	&&	\Br(\X_U)	\\
	&&	@VVV	\\
	\HH^1_\et(U,\PIC^0_{\X_U/U}) @>>> \HH^1_\et(U,\PIC_{\X_U/U}) @>>> 0.
\end{CD}
\end{equation}
The vertical map is the homomorphism
	$$\ker \left[ \HH^2_\et(\X_U,\G_m) \rightarrow
			\HH^0_\et(U,\RR^2\pi_\ast \G_m) \right]
	\longrightarrow \HH^1_\et(U,\RR^1\pi_\ast \G_m)$$
from the Leray spectral sequence
for $\pi: \X_U \rightarrow U$, since
the sheaf $\RR^1\pi_\ast \G_m$ can be identified with $\PIC_{\X_U/U}$
by~\cite[p.~203]{blr},
and $\RR^2\pi_\ast \G_m=0$ by Corollaire~3.2 in~\cite{grothbrauer}.

We restrict $\alpha'$ to an element of $\Br(\X_U)$,
map it downwards in~(\ref{bigLeray})
to obtain $\beta \in \HH^1_\et(U,\PIC_{\X_U/U})$,
and lift leftwards to obtain $\gamma \in \HH^1_\et(U,\PIC^0_{\X_U/U})$.
By~\cite[III.4.7]{milneetale},
$\gamma$ is the class of some torsor $\Y_U$
for $\PIC^0_{\X_U/U}=\E_U$ over $U$.

\begin{prop}
With notation as above,
assuming that there is a Brauer-Manin obstruction for $V$,
the class of $\Y_t$ is in $\Sha(\E_t)$ for all $t \in U(k)$,
and the Cassels-Tate pairing $\langle \X_t,\Y_t \rangle \in \Q/\Z$
is a nonzero constant independent of $t$.
\end{prop}

\begin{proof}
Suppose $t \in U(k)$.
If $\epsilon$ is a cohomology class of an \'etale sheaf on $U$ (or $\X_U$),
let $\epsilon_t$ denote the restriction to $t$ (resp.\ $\X_t$).
We have a diagram
\begin{equation}
\label{smallLeray}
\begin{CD}
	&&	\Br(\X_t)	\\
	&&	@VVV		\\
	\HH^1(k,\Pic^0(\Xbar_t)) @>>> \HH^1(k,\Pic(\Xbar_t)) @>>> 0.
\end{CD}
\end{equation}
analogous to~(\ref{bigLeray}),
in which the restrictions $\alpha'_t$ on the top
and $\gamma_t$ on the left map to $\beta_t$.
The class of the fiber $\Y_t$ equals $\gamma_t$.

Then $\ev(\alpha'_t,x)=\ev(\alpha,g(x))$ is a nonzero constant
independent of $x \in \X_t(\Ak)$ (and also independent of $t$).
The former is a sum of local functions $\ev_v(\alpha'_t,x_v)$
of independent arguments, so each of these must be constant.
Let $E=\X_t \times_k k_v$, which is an elliptic curve over $k_v$,
since $\X_t(k_v)\not=\emptyset$.
Let $\Ebar=E \times_{k_v} \kbar_v$.
A point in $\Pic^0(E)(k_v)$ is represented by
a difference of two points in $E(k_v)$,
and the pairing $\Br(E) \times \Pic(E) \rightarrow \Br(k_v)$
induced by evaluation
is compatible with the pairing
$\HH^1(k_v,\Pic \Ebar) \times \Pic^0(E)
\rightarrow \Br(k_v) \hookrightarrow \Q/\Z$
~\cite{lichtenbaum},
so the image $\beta_E$ of $\beta_t$ in $\HH^1(k_v,\Pic \Ebar)$
is in the kernel of the latter pairing.
It follows from Tate local duality (or a result of Witt if $v$ is real)
that $\beta_E=0$.

In the diagram analogous to~(\ref{smallLeray}), but over $k_v$
instead of $k$, the bottom surjection is now an isomorphism,
since the first map in the exact sequence
	$$\HH^0(k_v,\Pic(\Ebar)) \rightarrow \Z \rightarrow
		\HH^1(k_v,\Pic^0(\Ebar)) \rightarrow
		\HH^1(k_v,\Pic(\Ebar)) \rightarrow 0$$
is surjective.
Hence the image of $\gamma_t$ in $\HH^1(k_v,\Pic^0(\Ebar))$,
which maps to $\beta_E=0$, must itself be zero.
Thus $\gamma_t \in \Sha(\E_t)$.
By Section~\ref{casselstate},
$\langle \X_t,\Y_t \rangle = \ev(\alpha'_t,x)$
for any $x \in \X_t(\Ak)$, and the right hand side is independent of $t$.
\end{proof}

\section*{Acknowledgements}

We thank Sir Peter Swinnerton-Dyer for some discussions about cubic surfaces,
in particular for a suggestion
which led to the proof of Lemma~\ref{oneelement}.
We also thank Qing Liu for directing us to~\cite{liu}.


\end{document}